\documentclass[pdflatex,sn-mathphys-num]{sn-jnl}
\usepackage{amsmath,amssymb,amsthm}
\newtheorem{theorem}{Theorem}
\newtheorem*{maintheorem}{Main Theorem}
\newtheorem{lemma}[theorem]{Lemma}

\begin{document}
	\title{A martingale proof of simultaneous convex and smooth renorming}
	\author[1]{\fnm{Zhiyao} \sur{Fang}}
	\email{fzy73@stu.aqnu.edu.cn}
	\author[]{\fnm{Tingting} \sur{Hu}}
	\email{y25060048@stu.aqnu.edu.cn}
	\author[]{\fnm{Mengmeng} \sur{Bao}}
	\email{y25060053@stu.aqnu.edu.cn}
	\author*[1]{\fnm{Qi} \sur{Liu}}
	\email{liuq67@aqnu.edu.cn}
	\affil[1]{%
		\orgdiv{School of Mathematics and Statistics},
		\orgname{Anqing Normal University},
		\orgaddress{\city{Anqing}, \postcode{246133}, \state{Anhui}, \country{China}}
	}
	\abstract{We give an affirmative answer to Problem 198 of Guirao, Montesinos, and Zizler [\emph{Open Problems in the Geometry and Analysis of Banach Spaces}, Springer, Cham (2016), p.~90], also formulated by Pisier [\emph{Martingales in Banach Spaces}, Cambridge University Press (2016), p.~461]. Given equivalent norms with convexity of power type $p$ and smoothness of power type $q$, respectively, we construct a single equivalent norm with both properties using finite dyadic martingales. It can be chosen arbitrarily close to the given smooth norm in multiplicative distortion and preserves all linear symmetries shared by the two initial norms. The case $p=q=2$ with group invariance answers Problem 3.6 of Cabello S\'anchez, Ferenczi, and Randrianantoanina [S\~ao Paulo J. Math. Sci. \textbf{16} (2022), 406--458].}
	\maketitle
	\newpage
	
	\section{Introduction and main result}
	\label{sec:introduction}
	
	Guirao, Montesinos, and Zizler state the question as follows
	\cite[p.~90]{GMZ}.
	
	\begin{quote}
		\textbf{Problem 198.} If a Banach space $X$ admits a uniformly convex
		norm with modulus of convexity of power type $p$ and admits a
		uniformly smooth norm with modulus of smoothness of power type $q$,
		does $X$ admit a norm that simultaneously shares both these properties?
	\end{quote}
	
	Pisier formulates the question in terms of isomorphisms
	\cite[Section~11.6, p.~461]{PisierBook}.
	
	\begin{quote}
		\textbf{Problem:} If $B$ is both isomorphic to a $p$-uniformly smooth
		space and isomorphic to a $q$-uniformly convex one, is $B$ isomorphic
		to a space that is both $q$-uniformly convex and $p$-uniformly smooth?
	\end{quote}
	
	The two formulations are equivalent, with the roles of $p$ and $q$
	interchanged. A Banach-space isomorphism $T\colon X\to Y$ induces
	the equivalent norm $x\mapsto\lVert Tx\rVert_Y$ on $X$, whose
	moduli of convexity and smoothness are those of $Y$. Conversely,
	the identity map from $X$ to $X$ equipped with any equivalent norm
	is an isomorphism. Throughout the paper we follow the notation
	in the statement of Problem~198, with $p$ for convexity and $q$
	for smoothness.
	
	Pisier's interpolation argument
	\cite[Corollary~11.45 and p.~461]{PisierBook}
	gives an equivalent norm with both power-type properties, with
	a larger convexity exponent and a smaller smoothness exponent.
	The Main Theorem below preserves both prescribed exponents.
	
	The case in which both exponents equal two follows from the Hilbert space characterization discussed in \cite[p.~436]{CFR}. The same survey discusses the general problem and records that an almost transitive equivalent norm gives an affirmative answer in the superreflexive setting.
	
	Pisier's martingale renorming \cite{Pisier} relates power-type convexity to martingale cotype. Wenzel's proof of \cite[Theorem~3.6]{Wenzel} constructs a convex renorming by minimizing the difference between a terminal moment and a sum of increment moments. We use this construction with the terminal moment measured in the given smooth norm. The coefficient of the increment term is chosen so that the terminal moments of all near minimizers admit a uniform bound. This makes it possible to estimate a second difference of the infimum by translating the same near-minimizing martingale in both directions.
	
	A related result is Theorem~9.1 of Cheng, He, Liu, and Yu
	\cite{ChengEtAl}. It states that for every superreflexive Banach space
	$X$ there exist $2\le p<\infty$ and $1<q\le2$ such that the convex
	bodies in $X$ that are both $p$-uniformly convex and $q$-uniformly
	smooth are dense in the space of all convex bodies in the Hausdorff
	metric. In their proof, the authors use Asplund averaging and discuss
	the powers of the starting primal and dual norms. Our main theorem
	begins with a prescribed pair $(p,q)$ and constructs an equivalent
	norm by a dyadic martingale argument.
	
	Clarkson introduced uniformly convex spaces \cite{Clarkson}. Day studied uniform convexity in quotient spaces and its relation to smoothness in dual spaces \cite{Day}. Lindenstrauss introduced the modulus of smoothness and studied its relation to the modulus of convexity of the dual space \cite{Lindenstrauss}. For a systematic treatment of smoothness and equivalent norms, see Deville, Godefroy, and Zizler \cite{DGZ}.
	
	For a norm $B$ on a nonzero Banach space $X$, write $S_B=\{x\in X:B(x)=1\}$. Its moduli of convexity and smoothness are
	\begin{align*}
		\delta_B(\varepsilon)
		&=\inf_{\substack{u,v\in S_B\\B(u-v)\ge\varepsilon}}
		\left(1-B\left(\frac{u+v}{2}\right)\right),
		&&0<\varepsilon\le 2,\\
		\rho_B(t)
		&=\sup_{u,v\in S_B}
		\left(\frac{B(u+tv)+B(u-tv)}2-1\right),
		&&t>0.
	\end{align*}
	The moduli of convexity and smoothness are used to study both duality and
	the asymptotic geometry of Banach spaces. For duality, Lindenstrauss's
	formulas \cite{Lindenstrauss,LindenstraussTzafriri} express the smoothness
	modulus of either the original or the dual norm in terms of the convexity
	modulus of the other. Let $B^*$ denote the dual norm
	on $X^*$ and set $\delta_B(0)=\delta_{B^*}(0)=0$. Then for every $\tau>0$,
	\begin{align*}
		\rho_{B^*}(\tau)
		&=\sup_{0\le\varepsilon\le2}
		\left\{\frac{\tau\varepsilon}{2}-\delta_B(\varepsilon)\right\},\\
		\rho_B(\tau)
		&=\sup_{0\le\varepsilon\le2}
		\left\{\frac{\tau\varepsilon}{2}-\delta_{B^*}(\varepsilon)\right\}.
	\end{align*}
	In particular, a convexity bound of power type $r>1$ for either norm
	yields a smoothness bound of conjugate power type $r/(r-1)$ for the other.

	Convexity of power type $p$ means that $\delta_B(\varepsilon)\ge c\varepsilon^p$ for some $c>0$ and all $0<\varepsilon\le2$. Smoothness of power type $q$ means that $\rho_B(t)\le Kt^q$ for some $K>0$ and all $t>0$.
	
	In real dimension at least two, the Day--Nordlander bounds imply that the admissible exponents satisfy $1<q\le2\le p<\infty$ \cite{Nordlander,PisierBook}. In real dimension one, all norms are proportional, so the two estimates hold simultaneously whenever the hypotheses of Problem~198 hold. We work in the range $1<q\le2\le p<\infty$ below, over either scalar field.
	
	\stepcounter{theorem}
	\begin{maintheorem}
		Let $X$ be a real or complex Banach space and let $1<q\le2\le p<\infty$. Suppose $X$ admits equivalent norms $A$ and $B$ with convexity of power type $p$ and smoothness of power type $q$, respectively. For each $0<\theta<1$, there is an equivalent norm $N_\theta$ and constants $c_\theta,C_\theta>0$ such that
		\[
		\delta_{N_\theta}(\varepsilon)\ge c_\theta \varepsilon^p
		\quad (0<\varepsilon\le 2), \qquad
		\rho_{N_\theta}(t)\le C_\theta t^q
		\quad (t>0).
		\]
		Moreover,
		\begin{equation}
			\label{eq:comparison}
			(1-\theta)^{1/p}B(x)\le N_\theta(x)\le B(x)
			\qquad (x\in X).
		\end{equation}
	\end{maintheorem}
	
	\section{Proof and group invariance}
	\label{sec:proof}
	
	We first establish two midpoint estimates for the $p$th powers of $A$ and $B$ to control the martingale increments and the change in terminal cost under translation, respectively.
	
	\begin{lemma}
		\label{lem:power-estimates}
		Let $1<q\le2\le p<\infty$. If $A$ has convexity of power type $p$, there is $\alpha>0$ such that
		\begin{equation}
			\label{eq:convex-power}
			\frac{A(z+h)^p+A(z-h)^p}{2}-A(z)^p\ge \alpha A(h)^p
			\qquad (z,h\in X).
		\end{equation}
		If $B$ has smoothness of power type $q$, there is $L>0$ such that
		\begin{equation}
			\label{eq:smooth-power}
			\frac{B(z+h)^p+B(z-h)^p}{2}-B(z)^p
			\le L\bigl(B(z)+B(h)\bigr)^{p-q}B(h)^q
			\qquad (z,h\in X),
		\end{equation}
		where the factor of exponent zero is interpreted as $1$ when $p=q$.
	\end{lemma}
	
	\begin{proof}
		Inequality \eqref{eq:convex-power} is the standard midpoint characterization of power-type uniform convexity \cite[Section~10.1]{PisierBook}. To prove \eqref{eq:smooth-power}, put $a=B(z)$ and $b=B(h)$. The assertion is immediate if $b=0$. If $b>a$, the triangle inequality bounds the left-hand side by $2^p b^p$. Suppose $0<b\le a$ and write $s_\pm=B(z\pm h)$. Since $|s_\pm-a|\le b$ and $s_\pm\in[0,2a]$, Taylor's formula gives
		\[
		s_\pm^p\le a^p+pa^{p-1}(s_\pm-a)+C_p a^{p-2}b^2.
		\]
		The average of the linear terms is at most $pa^p\rho_B(b/a)$, which is bounded by a constant times $a^{p-q}b^q$. The remainder is bounded by the same power since $b\le a$ and $q\le2$.
	\end{proof}
	
	\begin{proof}[Proof of Main Theorem]
		\textbf{Step 1. Constructing the new norm.}
		Scale $A$ so that $A\le B$ and choose $\alpha$ as in \eqref{eq:convex-power}. A finite dyadic martingale $M$ starting at $x$ is represented by a finite binary tree with root $x$. Each nonterminal node is the midpoint of its two successors, and at each step either successor is selected with probability $1/2$. Let $M_k$ be the value at level $k$ and put $d_k=M_k-M_{k-1}$. For each such tree of depth $n$, put
		\begin{equation}
			\label{eq:cost}
			J(M)=\mathbb E B(M_n)^p-\theta\alpha\sum_{k=1}^n\mathbb E A(d_k)^p,
			\qquad
			F(x)=\inf_{M_0=x}J(M).
		\end{equation}
		Trees of depth zero are allowed. Applying \eqref{eq:convex-power} conditionally at each node and summing over levels gives
		\[
		\alpha\sum_{k=1}^n\mathbb E A(d_k)^p
		\le\mathbb E A(M_n)^p-A(x)^p.
		\]
		The normalization $A\le B$ therefore gives
		\begin{equation}
			\label{eq:terminal}
			J(M)\ge (1-\theta)\mathbb E B(M_n)^p+\theta A(x)^p.
		\end{equation}
		Since $\mathbb E M_n=x$, Jensen's inequality and the depth-zero tree give
		\begin{equation}
			\label{eq:Fbounds}
			(1-\theta)B(x)^p\le F(x)\le B(x)^p.
		\end{equation}
		
		Fix $x,h\in X$. For each $\eta>0$, take trees rooted at $x+h$ and $x-h$ whose costs are within $\eta$ of their respective infima. Pad with zero increments to a common depth and attach them to a root at $x$. The cost of the resulting tree is the average of the two costs minus $\theta\alpha A(h)^p$. Letting $\eta$ decrease to zero yields
		\begin{equation}
			\label{eq:lower-F}
			\frac{F(x+h)+F(x-h)}{2}-F(x)\ge\theta\alpha A(h)^p.
		\end{equation}
		By \eqref{eq:Fbounds}, the midpoint convex function $F$ is locally bounded above, hence continuous and convex. Scaling the trees shows that $F(\lambda x)=|\lambda|^pF(x)$, so $N_\theta=F^{1/p}$ is the gauge of the convex balanced set $\{F\le1\}$. Inequality \eqref{eq:Fbounds} proves that it is an equivalent norm and gives \eqref{eq:comparison}.
		
		\textbf{Step 2. Preserving the convexity exponent.}
		For $N_\theta(u)=N_\theta(v)=1$, take $x=(u+v)/2$, $h=(u-v)/2$ in \eqref{eq:lower-F}. Equivalence gives $A\ge mN_\theta$ for some $m>0$. Thus
		\[
		1-N_\theta\left(\frac{u+v}{2}\right)^p
		\ge \theta\alpha (m/2)^p N_\theta(u-v)^p.
		\]
		Using $1-s^p\le p(1-s)$ for $0\le s\le1$ proves the asserted convexity bound.
		
		\textbf{Step 3. Preserving the smoothness exponent.}
		For $\eta>0$, choose a tree rooted at $x$ with $J(M)\le F(x)+\eta$ and write $Z=M_n$. Translate this same tree by $h$ and by $-h$ to obtain competitors for $F(x+h)$ and $F(x-h)$. Because both translations leave every increment unchanged, the increment terms cancel when we subtract $J(M)$ from the average of the translated costs. Using \eqref{eq:smooth-power}, we obtain
		\[
		\frac{F(x+h)+F(x-h)}2-F(x)
		\le \eta+L\,\mathbb E \bigl(B(Z)+B(h)\bigr)^{p-q}B(h)^q.
		\]
		For $p>q$, the monotonicity of moment norms and the triangle inequality in $L_p$ imply
		\[
		\mathbb E\bigl(B(Z)+B(h)\bigr)^{p-q}
		\le \left((\mathbb E B(Z)^p)^{1/p}+B(h)\right)^{p-q}.
		\]
		For $p=q$ the expectation is simply $1$. By \eqref{eq:terminal}, the terminal moment is at most $(F(x)+\eta)/(1-\theta)$. Letting $\eta$ decrease to zero and using $F(x)\le B(x)^p$, we obtain
		\begin{equation}
			\label{eq:upper-F}
			\frac{F(x+h)+F(x-h)}2-F(x)
			\le L_\theta \bigl(B(x)+B(h)\bigr)^{p-q} B(h)^q.
		\end{equation}
		
		Finally, take $N_\theta(x)=N_\theta(y)=1$ and $h=ty$ for $0<t\le1$. By \eqref{eq:comparison}, the right-hand side of \eqref{eq:upper-F} is at most $K_\theta t^q$. The power mean inequality and the concavity of $s^{1/p}$ give
		\[
		\frac{N_\theta(x+ty)+N_\theta(x-ty)}2
		\le(1+K_\theta t^q)^{1/p}\le1+\frac{K_\theta}{p}t^q.
		\]
		For $t\ge1$ the triangle inequality gives $\rho_{N_\theta}(t)\le t\le t^q$. This proves the smoothness bound with exponent $q$.
	\end{proof}
	
	\begin{theorem}
		\label{thm:invariant}
		Under the assumptions of Main Theorem, suppose a group $G$ acts on $X$ by linear automorphisms and preserves both $A$ and $B$. Then $N_\theta$ can be chosen $G$-invariant.
	\end{theorem}
	
	\begin{proof}
		For $g\in G$, applying $g$ to every node of a tree rooted at $x$ gives a tree rooted at $gx$. The operation is bijective, with inverse obtained by applying $g^{-1}$. Each tree retains its cost in \eqref{eq:cost}, because the terminal moment uses $B$ and the increments use $A$. Consequently $F(gx)=F(x)$.
	\end{proof}
	
	\noindent\textbf{Problem 3.6} \cite[p.~437]{CFR}.
	Let $G$ be a bounded group of automorphisms of a Hilbert space $H$.
	Suppose $H$ admits a $G$-invariant equivalent norm with convexity
	modulus of power type $2$ and a $G$-invariant equivalent norm with
	smoothness modulus of power type $2$. Does $H$ admit a single
	$G$-invariant equivalent norm whose two moduli are both of power type $2$?
	
	Theorem~\ref{thm:invariant} answers Problem~3.6 affirmatively. Take $A$
	and $B$ to be the two $G$-invariant norms in the question, respectively,
	and apply the theorem with $X=H$ and $p=q=2$. The resulting norm
	$N_\theta$ is $G$-invariant and satisfies
	$\delta_{N_\theta}(\varepsilon)\ge c_\theta\varepsilon^2$ and
	$\rho_{N_\theta}(t)\le C_\theta t^2$.

	\section*{Acknowledgements}
	
	The authors thank the Functional Analysis Research Team at the School of Math-
	ematics and Statistics of Anqing Normal University for discussing the manuscript
	and helping to resolve difficulties in the proof.
	
	\section*{Declarations}
	
	\paragraph{Competing interests} The authors declare that they have no competing interests.
	
	\paragraph{Funding} This research received no specific funding.
	
	\paragraph{Data availability} No datasets were generated or analysed for this study.


\begin{thebibliography}{99}
		
		\bibitem{CFR}
		Cabello S\'anchez, F., Ferenczi, V., Randrianantoanina, B.:
		On Mazur rotations problem and its multidimensional versions.
		S\~ao Paulo J. Math. Sci. \textbf{16}, 406--458 (2022)
		
		\bibitem{ChengEtAl}
		Cheng, L., He, W., Liu, C., Yu, Z.:
		On strengthened versions of Klee's convex body problem in Banach spaces.
		Preprint, arXiv:2606.07370v1 (2026)
		
		\bibitem{Clarkson}
		Clarkson, J.A.:
		Uniformly convex spaces.
		Trans. Amer. Math. Soc. \textbf{40}, 396--414 (1936)
		
		\bibitem{Day}
		Day, M.M.:
		Uniform convexity in factor and conjugate spaces.
		Ann. of Math. (2) \textbf{45}, 375--385 (1944)
		
		\bibitem{DGZ}
		Deville, R., Godefroy, G., Zizler, V.:
		Smoothness and Renormings in Banach Spaces.
		Pitman Monographs and Surveys in Pure and Applied Mathematics, vol.~64.
		Longman Scientific \& Technical, Harlow (1993)
		
		
		\bibitem{GMZ}
		Guirao, A.J., Montesinos, V., Zizler, V.:
		Open Problems in the Geometry and Analysis of Banach Spaces.
		Springer, Cham (2016)
		
		\bibitem{Lindenstrauss}
		Lindenstrauss, J.:
		On the modulus of smoothness and divergent series in Banach spaces.
		Michigan Math. J. \textbf{10}, 241--252 (1963)
		
		\bibitem{LindenstraussTzafriri}
		Lindenstrauss, J., Tzafriri, L.:
		Classical Banach Spaces II. Function Spaces.
		Ergebnisse der Mathematik und ihrer Grenzgebiete, vol.~97.
		Springer, Berlin (1979)
		
		\bibitem{Nordlander}
		Nordlander, G.:
		The modulus of convexity in normed linear spaces.
		Ark. Mat. \textbf{4}, 15--17 (1960)
		

		
		\bibitem{Pisier}
		Pisier, G.:
		Martingales with values in uniformly convex spaces.
		Israel J. Math. \textbf{20}, 326--350 (1975)
		
		\bibitem{PisierBook}
		Pisier, G.:
		Martingales in Banach Spaces.
		Cambridge Studies in Advanced Mathematics, vol.~155.
		Cambridge University Press, Cambridge (2016)
		
		\bibitem{Wenzel}
		Wenzel, J.:
		Strong martingale type and uniform smoothness.
		J. Convex Anal. \textbf{12}, 159--171 (2005)
		
	\end{thebibliography}
\end{document}